\documentclass[sn-mathphys,Numbered]{sn-jnl}

\usepackage{graphicx}%
\usepackage{multirow}%
\usepackage{amsmath,amssymb,amsfonts}%
\usepackage{amsthm}%
\usepackage{mathrsfs}%
\usepackage[title]{appendix}%
\usepackage{xcolor}%
\usepackage{textcomp}%
\usepackage{manyfoot}%
\usepackage{booktabs}%
\usepackage{algorithm}%
\usepackage{algorithmicx}%
\usepackage{algpseudocode}%
\usepackage{listings}%
\theoremstyle{thmstyleone}%
\newtheorem{theorem}{Theorem}

\newtheorem{proposition}[theorem]{Proposition}
\newtheorem{lemma}[theorem]{Lemma}
\newtheorem{corollary}[theorem]{Corollary}

\theoremstyle{thmstyletwo}%
\newtheorem{remark}{Remark}%

\theoremstyle{thmstylethree}%

\def\bR{\mathbb{R}}

\def\bN{\mathbb{N}}

\def\ee{\mathbb{E}}
\def\prob{\mathbb{P}}

\begin{document}

\title{Long run stochastic control problems with general discounting}

\author{\fnm{Łukasz} \sur{Stettner}} {\email{stettner@impan.pl}}

\affil{\orgdiv{Institute of Mathematics}, \orgname{Polish Acad. Sci.}, \orgaddress{\street{Sniadeckich 8}, \city{Warsaw}, \postcode{00-656},  \country{Poland}}}

\abstract
{Controlled discrete time Markov processes are studied first with long run general discounting functional. It is shown that optimal strategies for average reward per unit time problem are also optimal for average generally discounting functional. Then long run risk sensitive reward functional with general discounting is considered. When risk factor is positive then optimal value of such reward functional is dominated by the reward functional  corresponding to the long run risk sensitive control. In the case of negative risk factor we get an asymptotical result, which says that optimal average reward per unit time control is nearly optimal for long run risk sensitive reward functional with general discounting, assuming that risk factor is close to $0$. For this purpose we show in Appendix upper estimates for large deviations of weighted empirical measures, which are of independent interest.}



\maketitle

\section{Introduction}\label{intro}

Let $(X_n)$ be a discrete Markov process on $(\Omega, F,(F_n))$ taking values in a locally compact separable metric space $E$ endowed with Borel $\sigma$ field ${\cal E}$. The process has transition operator (kernel) $\prob^a(x,\cdot)$, where $x\in E$ and $a\in U$ a fixed compact set of control parameters. We shall assume that the mapping $E\times U\ni (x,a)\mapsto \prob^a(x,D)$ is Borel measurable for $D\in {\cal E}$.  Markov process is controlled using sequence $V=(a_0,a_1,\ldots, a_n\ldots)$ such that $a_n\in U$ is $F_n$ adapted. We denote by $\ee_{k,x}^V$ expected value corresponding to controlled Markov process with the use of sequence $V$ and starting at time $k$ from $x$. When $k=0$ we simply write $\ee_{x}^V$.  We want to maximize the following long run average generally discounted functional
\begin{equation} \label{fun1}
J_{k,x}^g(V)=\liminf_{n\to \infty}\left(\sum_{i=k}^{n+k-1}\phi(i)\right)^{-1} \ee_{k,x}^V\left[\sum_{i=k}^{n+k-1}\phi(i)c(X_{i-k},a_{i-k})\right],
\end{equation}
where $\phi: \bN\to [0,1]$ is a generalized discount factor such that: $\phi(0)=1$, $n\to \phi(n)$ is nonincreasing and $\phi(n+k)\geq\phi(n)\phi(k)$, for $n,k\in \bN$.
We shall assume furthermore that $\sum_{k=0}^\infty \phi(k)=\infty$.
Generally we have in mind generalized discount factors such that $n\to \phi(n)$ is strictly decreasing and $\phi(n+k)>\phi(n)\phi(k)$, for $n,k\in \bN$.
Typical example for such discount rate is a generalized hyperbolic discount factor of the form $\phi(i)=(1+hi)^{-{r\over h}}$, where $h>0, r>0$ and ${r\over h}\leq 1$. It plays an important role in economics see e.g. \cite{Loew1992}, \cite{Prelec2004}, \cite{Rohde2010}, \cite{Al2006} and references therein. The problems with non-exponential discounting are time inconsistent (see sections 7 and 17 of \cite{Bjork2021} for examples of such problems or \cite{Bau2021}). General hyperbolic discounting was considered in the case of repeated games in \cite{Obara2017}. Finite horizon consumption and portfolio decisions with stochastic
hyperbolic discounting was studied in \cite{Zoua}. The necessity to use hyperbolic discounting working with market date was shown in \cite{Salois}. It appears that working with far distant futures we have to consider hyperbolic discounting as was noted in the paper \cite{Anch}. Further motivations to study hyperbolic discount one can find in \cite{Pen}. In the paper we shall study long time behaviour of functional \eqref{fun1} and the formulated below functional \eqref{fun2} with general discounting. We shall compare with the case $\phi\equiv 1$, which corresponds to the undiscounted problems.
Moreover we assume that the mapping $E\times U\ni (x,a) \mapsto c(x,a)$ is continuous and bounded. We are also interested in maximization of the following long run risk sensitive generally discounted functional
\begin{equation} \label{fun2}
I_{k,x}^{\gamma,g}(V)=\liminf_{n\to \infty}\left(\gamma\sum_{i=k}^{n+k-1}\phi(i)\right)^{-1}  \ln\left( \ee_{k,x}^V\left[\exp\left\{\gamma\sum_{i=k}^{n+k-1}\phi(i)c(X_{i-k},a_{i-k})\right\}\right]\right)
\end{equation}
with $\gamma \neq 0$. This functional measures not only averaged discounted reward but also its higher moments with suitable weights, in particular variance which is considered as a measure of risk.
The problems with functionals \eqref{fun1} and in particular with functional \eqref{fun2} are difficult to study and are time inconsistent. The purpose of the paper is to show that under certain assumptions, usually for the computational purposes, when we  make suitable discretizations of the problems, the problems can be fully or partially solved using undiscouted versions of the functionals \eqref{fun1} and \eqref{fun2}.
 In the undiscounted case, when $\phi\equiv 1$ instead of \eqref{fun1} we maximize average reward per unit time functional
\begin{equation} \label{fun3}
J_{x}(V)=\liminf_{n\to \infty}{1\over n} \ee_x^V\left[\sum_{i=0}^{n-1}c(X_{i},a_{i})\right],
\end{equation}
while in the case of \eqref{fun2} we maximize long run risk sensitive functional
\begin{equation} \label{fun4}
I_{x}^{\gamma}(V)=\liminf_{n\to \infty}{1\over n\gamma}  \ln\left( \ee_x^V\left[\exp\left\{\gamma\sum_{i=0}^{n-1}c(X_{i},a_{i})\right\}\right]\right).
\end{equation}
We shall study reward functionals \eqref{fun1}-\eqref{fun2} comparing them to \eqref{fun3}-\eqref{fun4} respectively.  Our fundamental assumption will be so called uniform ergodicity condition
\begin{equation}\label{A.1} \tag{A.1}
\sup_{x,x'\in E} \sup_{a,a'\in U} \sup_{D \in {\cal E}} \prob^a(x,D)-\prob^{a'}(x',D)=\Delta<1.
\end{equation}
The assumption \eqref{A.1} is rather restrictive in the case of locally compact state spaces but is commonly satisfied for nondegenerate processes on compact state spaces. Given stationary Markov control i.e. Borel measurable mapping $u:E\to U$, Markov process with transition operator $\prob^{u(x)}(x,\cdot)$ is uniformly ergodic and by section 5.5 of \cite{Doob} it has a unique invariant measure $\mu_u \in{\cal P}(E)$, where ${\cal P}(E)$ denotes the space of probability measures on $E$.
Another important assumption is a weak continuity of the controlled transition kernels (sometimes called Feller property of the controlled kernels)   that is
\begin{equation} \label{A.2} \tag{A.2}
the\ \  mapping \ \  E\times U\ni (x,a)\mapsto \prob^af(x) \ \ is \ \ continuous,
\end{equation}
where $\prob^af(x)=\int_E f(y)\prob^a(x,dy)$ for $f\in C(E)$ - the space of continuous bounded functions on $E$.  Under assumption \eqref{A.1} and \eqref{A.2} the functional \eqref{fun3} was studied  in particular in \cite{HerLas1996} (see also the references therein). The functional \eqref{fun4} was investigated in \cite{DiMSte1999} and in \cite{{DiMSte2000}}. Under the above assumptions bounded solutions to so called Bellman equations corresponding to such functionals were obtained. Boundedness of solutions to the Bellman equations plays an important role in the analysis of problems with general discount factor. We are aware of the fact that  \eqref{fun4} was recently studied in a number of papers under more general assumptions: see \cite{Jas2007} using Bellman inequalities, \cite{PitSte2016} for portfolio applications, \cite{WeiChen2021} for risk sensitive nonzero sum games,  \cite{GolPal2022} for discounted continuous time problem, \cite{GuoH2021} for continuous time with unbounded reward functions or for countable state spaces in \cite{BisPra2022} and references therein. In all these cases solutions to suitable Bellman equations were frequently unbounded (see also \cite{Stettner} with suitable examples and comments). We don't know results concerning the problems with functionals \eqref{fun1} or \eqref{fun2}. These problems are time inconsistent (see \cite{Bjork2021} for detailed study of time inconsistent problems). One of potential approaches to such problems is based on extension of the state space by adding time variable. We first follow this approach studying suitable Bellman equations for controlled state - time processes. Then we use the results on undiscounted problems to study discounted functionals \eqref{fun1} and \eqref{fun2}. We show that optimal controls for \eqref{fun3} are also optimal for \eqref{fun1}. The case of \eqref{fun2} is more complicated. We show that when risk factor is positive the optimal value of the functional \eqref{fun2} is dominated by that of \eqref{fun4}. When risk factor is negative we get an asymptotical result. Namely optimal control for \eqref{fun3} is nearly optimal for \eqref{fun2} provided that risk factor is close to $0$. This result is based on large deviation upper bounds for weighted empirical measures. We formulate and show such result in appendix extending there old papers \cite{DV1} and \cite{DV3}, which is of independent interest. To obtain this result we have to assume \eqref{A.3}, which is even stronger than \eqref{A.1} and \eqref{A.2}, but is necessary to obtain pointed above large deviation result, when the state space is not compact. Such assumptions are commonly acceptable for models with compact state spaces, in particular when we want to make computations.

\section{Long run average discounted functional problem.} \label{Average}

Let
\begin{equation}\label{v1}
v^g(k,x):= \sup_V J_{k,x}^g(V)
\end{equation}
be the value function of the reward functional \eqref{fun1}. We have the following Bellman equation corresponding to this functional
\begin{equation}\label{B1}
w(i,x)=\sup_{a\in U} \left[\phi(i)(c(x,a)-\lambda(i)) + \int_E w(i+1,y)\prob^a(x,dy)\right]
\end{equation}
for $i\geq k$.
We are looking for a sequence $(\lambda(i))$ and continuous bounded function $w$ for which \eqref{B1} is satisfied. For this purpose we introduce the span norm $\|f\|_{sp}:=\sup_{x\in E} f(x) -\inf_{x'\in E} f(x')$, for $f\in C(E)$ and the quotient space $C_{sp}=C(\bN\times E)\slash \sim$, where for $f,g \in C(\bN\times E)$ we have $f\sim g \equiv \forall_{i\in \bN} \ \|f(i,\cdot)-g(i,\cdot)\|_{sp}=0$. We endow the space $C_{sp}$ with the norm
$\mid\|f\|\mid := sup_{i\in \bN} \|f(i,\dot)-g(i,\cdot)\|_{sp}$.
One can show that $C_{sp}$ with the norm $\mid\|\cdot\|\mid$ is a Banach space. We have

\begin{theorem}\label{thm1}
Under \eqref{A.1} and \eqref{A.2} there are solutions $w\in C_{sp}$ and sequence $(\lambda(i))$ to the equation \eqref{B1}. Moreover
\begin{equation}\label{vv1}
v^g(k,x)=\liminf_{n\to \infty} {\sum_{i=k}^{n+k-1} \lambda(i) \phi(i) \over \sum_{i=k}^{n+k-1} \phi(i)}
\end{equation}
and the above limit does not depend on $k\in \bN$ and $x\in E$. The optimal strategy $\hat{V}$ for the reward functional $J_{k,x}^g$ is of the form $\hat{V}=(\hat{u}_k(X_0),\hat{u}_{k+1}(X_1),\ldots \hat{u}_{k+i}(X_i),\ldots)$, where $\hat{u}_i$, for $i\geq k$ is a Borel measurable selector for which equality on the right hand side of \eqref{B1} is attained.
\end{theorem}

\proof
Consider the transformation $\Phi$ of the space $C_{sp}$ defined as follows for $i\geq k$, $w\in C_{sp}$
\begin{equation}\label{tr}
\Phi w(i,x):=\sup_{a\in U} \left[\phi(i)c(x,a)+ \int_E w(i+1,y) \prob^a(x,dy)\right].
\end{equation}
Let $w_1, w_2 \in C_{sp}$. For  $x,x'\in E$, $i\geq k$, by \eqref{A.2} there are $a_1,a_2 \in U$ such that
\begin{equation}\label{r1}
\Phi w_1(i,x)=\phi(i)c(x,a_1)+ \int_E w_1(i+1,y) \prob^{a_1}(x,dy)
\end{equation}
\begin{equation}\label{r2}
\Phi w_2(i,x')=\phi(i)c(x',a_2)+ \int_E w_2(i+1,y) \prob^{a_2}(x',dy).
\end{equation}
Then by \eqref{r1} and \eqref{r2} we have
\begin{eqnarray}\label{r3}
&& \Phi w_1(i,x)-\Phi w_2(i,x)- (\Phi w_1(i,x')-\Phi w_2(i,x'))\leq  \nonumber \\
&&\int_E (w_1(i+1,y)-w_2(i+1,y)) \left(\prob^{a_1}(x,dy)-\prob^{a_2}(x',dy)\right)\leq \nonumber \\
&&\sup_{y\in E} (w_1(i+1,y)-w_2(i+1,y))\left(\prob^{a_1}(x,H)-\prob^{a_2}(x',H)\right)- \nonumber \\
&&\inf_{y\in E} (w_1(i+1,y)-w_2(i+1,y))\left(\prob^{a_1}(x,H^c)-\prob^{a_2}(x',H^c)\right)=\nonumber \\
&& \|w_1(i+1,\cdot)-w_2(i+1,\cdot)\|_{sp} \left(\prob^{a_1}(x,H)-\prob^{a_2}(x',H)\right)\leq \nonumber \\
&& \|w_1(i+1,\cdot)-w_2(i+1,\cdot)\|_{sp} \Delta,
\end{eqnarray}
where $H$ comes from the Hahn decomposition of $\prob^{a_1}(x,\cdot)-\prob^{a_2}(x',\cdot)$.
Therefore we have that
\begin{equation}\label{r4}
\mid \|\Phi w_1-\Phi w_2\|\mid \leq \Delta \mid \|w_1-w_2\|\mid
\end{equation}
and there is a fixed point $w\in C_{sp}$ such that $\|w-\Phi w\|_{sp}=0$, which means that there is a sequence $(\lambda(i))$ such that  $w(i,x)+\lambda(i)\phi(i)=\Phi w(i,x)$ for $x\in E$ and $i\geq k$. Consequently $w$ and $(\lambda(i))$ form a solution to \eqref{B1}. Fix $x_0\in E$ and define $\bar{w}(i,x):=w(i,x)-w(i,x_0)$, for $i\geq k$ and $x\in E$. Since $w\in C_{sp}$ we therefore have that $\sup_{i \geq k}\sup_{x\in E} |\bar{w}(i,x)|\leq \mid \|w\|\mid<\infty$. Consequently for $\tilde{\lambda}(i)\phi(i):=\lambda(i)\phi(i)+w(i,x_0)-w(i+1,x_0)$ we have
\begin{equation}\label{r5}
\bar{w}(i,x)+\tilde{\lambda}(i)\phi(i)=\sup_{a\in U} \left[\phi(i)c(x,a)+ \int_E \bar{w}(i+1,y)\prob^a(x,dy)\right]
\end{equation}
and for any control $V$
\begin{equation}\label{r6}
\bar{w}(k,x)+ \sum_{i=k}^{n+k-1} \tilde{\lambda}(i)\phi(i)\geq \ee_x^V\left[\sum_{i=k}^{n+k-1}\phi(i)c(X_{i-k},a_{i-k})+\bar{w}(n+k, X_n)\right],
\end{equation}
with equality when $a_{i-k}=u_i(X_{i-k})$, where $u_i$ is a Borel measurable selector of the right hand side of \eqref{r5}. Dividing \eqref{r6} by $\sum_{i=k}^{n+k-1}\phi(i)$ and letting $n\to \infty$ we obtain \eqref{vv1}. Notice that the limit does not depend on initial values of the sequences $(\tilde{\lambda}(i))$ and
$(\phi(i))$, which completes the proof.
\endproof

\begin{remark}\label{rem1}
It is clear from the proof that optimal strategies for the reward functional \eqref{fun1} are time dependent. The problem in a stationary (time independent) form is therefore time inconsistent (see more about such problems in \cite{Bjork2021}).  To make it time consistent it we had to add time variable. In what follows we shall prove that optimal strategy for \eqref{fun1} can be chosen also in a stationary (time independent) form.
\end{remark}
In the case of functional \eqref{fun3} by analogy to the first part of the proof of Theorem \ref{thm1} (also see \cite{HerLas1996})  we characterize of the value function
\begin{equation}\label{v3}
v(x):=\sup_V J_x(V),
\end{equation}
using  Bellman equation. We are looking for a function $w\in C(E)$ and a constant $\lambda$ such that
\begin{equation}\label{B3}
w(x)=\sup_{a\in U} \left[c(x,a)-\lambda + \int_E w(y)\prob^a(x,dy)\right]
\end{equation}
holds for $x\in E$. We immediately have
\begin{corollary}\label{cor1}
Under \eqref{A.1} and \eqref{A.2} there are $w \in C(E)$ and a constant $\lambda$ for which equation \eqref{B3} is satisfied. Furthermore $\|w\|_{sp}\leq {1 \over 1-\Delta}$,
\begin{equation}\label{vv3}
v(x)=\lambda
\end{equation}
and the control $\hat{V}=(\hat{u}(X_0),\hat{u}(X_1),\ldots,\hat{u}(X_n),\ldots)$, where $\hat{u}$ is a Borel selector to the supremum on the right hand side of \eqref{B3}, is optimal for the functional \eqref{fun3}.
\end{corollary}
\proof The proof in the case $\phi\equiv 1$  is much simpler than that of Theorem \ref{thm1}. Namely, we
 consider the transformation $\Phi$ of C(E) of the form
\begin{equation}\label{tr3}
\Phi w(x):=\sup_{a\in U} \left[c(x,a)+ \int_E w(y) \prob^a(x,dy)\right].
\end{equation}
By analogy to the first part of the proof of Theorem \ref{thm1} for $w_1,w_2\in C(E)$ we have
\begin{equation}
\|\Phi w_1-\Phi w_2\|_{sp}\leq \Delta \|w_1-w_2\|_{sp}.
\end{equation}
Therefore there is $w\in C(E)$ such that $\|w-\Phi w\|_{sp}=0$ and consequently for $w$ and $\lambda=w(x)-\Phi w(x)$, for any $x\in E$, the equation \eqref{B3} is satisfied. Moreover for $\Phi^0 0\equiv 0$ we have
\begin{equation}\label{r35'}
\|w\|_{sp}=\lim_{n\to \infty} \|\sum_{i=1}^{n}(\Phi^i 0- \Phi^{i-1}0)\|_{sp}\leq {1 \over 1-\Delta}.
\end{equation}
Furthermore for any strategy $V$ iterating \eqref{B3} we obtain
\begin{equation}\label{r36}
{w}(x)+ n\lambda \geq \ee_x^V\left[\sum_{i=0}^{n-1}c(X_{i},a_{i})+{w}(X_n)\right]
\end{equation}
 with equality for the strategy $\hat{V}$. Therefore
\eqref{vv3} follows. Optimality of the strategy $\hat{V}$ we obtain in a standard way (see \cite{HerLas1996}).
\endproof

\begin{remark}\label{rem2}
 Notice that since $w$ in \eqref{r36} is bounded we have even more, namely
\begin{equation}
\lambda=\sup_V \limsup_{n\to \infty} {1\over n} \ee_x^V\left[\sum_{i=0}^{n-1}c(X_{i},a_{i})\right].
\end{equation}
\end{remark}
We have now the following rather unexpected result

\begin{theorem}\label{thm2}
Under \eqref{A.1} and \eqref{A.2} the strategy $\hat{V}$ defined in Corollary \ref{cor1} is optimal for \eqref{fun1} and for $k\in \bN$ and $x\in E$ we have
$v^g(k,x)=\lambda$.
\end{theorem}

\proof
It is clear that $v^g(k,x)\geq J_{k,x}^g(\hat{V})$, for $k\in \bN$ and $x\in E$. Consider solutions $w\in C(E)$ and $\lambda$ to the Bellman equation \eqref{B3}. Without loss of generality we may assume that $w(x)\geq 0$ for $x \in E$ (otherwise we may consider $\bar{w}(x)= w(x)-\inf_{y\in E} w(y)$ which is also a solution to \eqref{B3} with the same constant $\lambda$). For any strategy $V=(a_0,a_1,\ldots, a_n\ldots)$ from \eqref{B3} we have
\begin{equation}
c(X_i,a_i)\leq w(X_i) + \lambda - \int_E w(y) \prob^{a_i}(X_i,dy).
\end{equation}
Then taking into account that $w(y)\geq 0$ for $y\in E$ and $\phi(i)$ is nonincreasing we have
\begin{eqnarray}\label{eqq1}
&&\ee_{k,x}^V\left[\sum_{i=k}^{n+k-1} \phi(i) c(X_{i-k},a_{i-k})\right]\leq \nonumber \\
&& \ee_{k,x}^V\left[\lambda \sum_{i=k}^{n+k-1} \phi(i) + \sum_{i=k}^{n+k-1} (w(X_{i-k})- \int_E w(y) \prob^{a_{i-k}}(X_{i-k},dy))\phi(i)\right]= \nonumber \\
&& \lambda \sum_{i=k}^{n+k-1} \phi(i)+ \sum_{i=k}^{n+k-1} \phi(i)\ee_{k,x}^V\left[w(X_{i-k})-w(X_{i-k+1})\right]=\nonumber \\
&& \lambda \sum_{i=k}^{n+k-1} \phi(i)+ \phi(k)w(x)- \ee_{k,x}^V\left[\sum_{i=k+1}^{n+k-1}w(X_{i-k})(\phi(i-1)-\phi(i))+ \right. \nonumber \\
&& \left. \phi(n+k-1)w(X_{n})\right]\leq  \lambda \sum_{i=k}^{n+k-1} \phi(i) + \phi(k)w(x).
\end{eqnarray}
Therefore $J_{k,x}^g({V})\leq \lambda$. On the other hand we have that
\begin{equation}\label{eqq2}
0\leq \ee_{k,x}^V\left[\sum_{i=k+1}^{n+k-1}w(X_{i-k})(\phi(i-1)-\phi(i))\right]\leq \|w\|(\phi(k)-\phi(n+k-1))\leq \|w\|,
\end{equation}
so that dividing \eqref{eqq1} by $\sum_{i=k}^{n+k-1} \phi(i)$, taking into account \eqref{eqq2} and letting $n\to \infty$ for strategy $\hat{V}$ we obtain  $J_{k,x}^g(\hat{V})=\lambda$. This completes the proof.
\endproof
\begin{remark}\label{rem3}
 Since $w$ is bounded we have shown in the proof above  even more, namely we have
\begin{equation}
\lambda=\sup_V\limsup_{n\to \infty}\left(\sum_{i=k}^{n+k-1}\phi(i)\right)^{-1} \ee_x^V\left[\sum_{i=k}^{n+k-1}\phi(i)c(X_{i-k},a_{i-k})\right].
\end{equation}
\end{remark}
We shall now study solutions to so called additive Poisson equations corresponding to Markov controls $V_m^k=(u_k(X_0),u_{k+1}(X_1),\ldots,u_{k+i}(X_i),\ldots)$ where $u_i:E\mapsto U$ is Borel measurable. Consider the space $B(E)$ consisting of bounded Borel measurable functions on $E$  and the quotient space $B_{sp}=B(\bN\times E)\slash \sim$, where for $f,g \in B(\bN\times E)$ we have $f\sim g \equiv \forall_{i\in \bN} \  \|f(i,\cdot)-g(i,\cdot)\|_{sp}=0$. We endow the space $B_{sp}$ with the norm
$\mid\|f\|\mid := sup_{i\in \bN} \|f(i,\dot)-g(i,\cdot)\|_{sp}$.
One can show that $B_{sp}$ with the norm $\mid\|\cdot\|\mid$ is also a Banach space.
Having fixed Markov control $V_m^k$ we are looking for  solutions $w^{V_m^k}$ and $(\lambda^{V_m^k}(i))$ to the equation
\begin{equation}\label{P1}
w^{V_m^k}(i,x)=\left[\phi(i)(c(x,u_{i}(x)-\lambda^{V_m^k}(i)) + \int_E w^{V_m^k}(i+1,y)\prob^{u_{i}(x)}(x,dy)\right]
\end{equation}
for $i\geq k$.
We have
\begin{corollary}\label{cor2}
Under \eqref{A.1} for a given Markov control $V_m^k$ there are solutions $w^{V_m^k}\in B_{sp}$ and sequence $(\lambda^{V_m^k}(i))$ to the equation \eqref{P1}. Moreover
\begin{equation}\label{vv11}
J^g_{k,x}(V_m^k)=\liminf_{n\to \infty} {\sum_{i=k}^{n+k-1} \lambda^{V_m^k}(i) \phi(i) \over \sum_{i=k}^{n+k-1} \phi(i)}
\end{equation}
and the above limit does not depend on $k\in \bN$ and $x\in E$.
\end{corollary}
\proof We follow the proof of Theorem \ref{thm1}. Namely we define the operator
\begin{equation}\label{trp}
\Phi w(i,x):= \left[\phi(i)c(x,u_i(x))+ \int_E w(i+1,y) \prob^{u_i(x)}(x,dy)\right]
\end{equation}
on the space $B_{sp}$ and using \eqref{A.1} we obtain contraction property, from which as in the proof of Theorem \ref{thm1} we have existence of solutions to \eqref{P1}.
\endproof
In particular case, when $V_m=(u(X_0),u(X_1),\ldots, u(X_i),\ldots)$, where $u:E\to U$ is Borel measurable we have a stationary Markov control. To simplify notation we shall identify control $V_m$ with function $u$ and as in Corollary \ref{cor2} we obtain
\begin{corollary}\label{cor3}
Under \eqref{A.1} for a stationary Markov control $u$ there is a bounded function $w^u\in B(E)$ and a constant $\lambda^u$ such that the Poisson equation
\begin{equation}\label{P3}
w^{u}(x)=\left[(c(x,u(x))-\lambda^u) + \int_E w^u(y)\prob^{u(x)}(x,dy)\right]
\end{equation}
is satisfied. Moreover $\|w^u\|_{sp}\leq {1 \over 1-\Delta}$ and
\begin{equation}\label{vvv1}
\lambda^u=J_x(u)=\lim_{n\to \infty}{1\over n} \ee_x^u\left[\sum_{i=0}^{n-1} c(X_i,u(X_i))\right].
\end{equation}
\end{corollary}
Next Corollary is a version of Theorem \ref{thm2} in the case of stationary Markov control $u$.
\begin{corollary}\label{cor4}
Under \eqref{A.1} for stationary Markov control $V_m=(u(X_0),u(X_1),\ldots,u(X_i),\ldots)=u$ we have
\begin{equation}\label{equale}
J_{k,x}^g(V_m)=\lambda^u.
\end{equation}
\end{corollary}
\proof
By analogy to \eqref{eqq1} from \eqref{P3} we have
\begin{equation}
c(X_i,u(X_i))= w^{u}(X_i)+\lambda^u - \int_E w^u(y)\prob^{u(X_i)}(X_i,dy)
\end{equation}
and therefore
\begin{eqnarray}\label{eqqn}
&&\ee_{k,x}^{V_m}\left[\sum_{i=k}^{n+k-1} \phi(i) c(X_{i-k},u(X_{i-k})\right]= \nonumber \\
&& \ee_{k,x}^{V_m}\left[\lambda^u \sum_{i=k}^{n+k-1} \phi(i) + \sum_{i=k}^{n+k-1} (w^u(X_{i-k})- \int_E w^u(y) \prob^{u(X_{i-k})}(X_{i-k},dy))\phi(i)\right]= \nonumber \\
&& \lambda^u \sum_{i=k}^{n+k-1} \phi(i)+ \sum_{i=k}^{n+k-1} \phi(i)\ee_{k,x}^{V_m}\left[w^u(X_{i-k})-w^u(X_{i-k+1})\right]=\nonumber \\
&& \lambda^u \sum_{i=k}^{n+k-1} \phi(i)+ \phi(k)w^u(x)- \ee_{k,x}^{V_m}\left[\sum_{i=k+1}^{n+k-1}w^u(X_{i-k})(\phi(i-1)-\phi(i))+ \right. \nonumber \\
&& \left. \phi(n+k-1)w^u(X_{n})
\right].
\end{eqnarray}
Without loss of generality we may assume in \eqref{P3}  that $w^u\geq 0$.
Then
\begin{equation}
0\leq \sum_{i=k+1}^{n+k-1}w^u(X_{i-k})(\phi(i-1)-\phi(i))\leq \|w^u\|
\end{equation}
and dividing \eqref{eqqn} by $\sum_{i=k}^{n+k-1} \phi(i)$ and letting $n\to \infty$ we obtain \eqref{equale}.
\endproof

We shall need also the following assumption:

there is $\eta\in {\cal P}(E)$ such that
\begin{equation}\label{A.3} \tag{A.3}
\prob^a(x,B)=\int_B p(x,a,y)\eta(dy) \ \ for \ \ B\in{\cal E} \ and  \ x\in E,
\end{equation}
the mapping $E\times U \ni (x,a) \to p(x,a,y)$ is continuous for $y\in E$ and there is $M>0$ such that ${1\over M} \leq p(x,a,y) \leq  M$ for $x,y\in E$ and $a\in U$.

\begin{remark}\label{rem4}
One can easily notice that under \eqref{A.3} both \eqref{A.1} and \eqref{A.2} are satisfied. In fact, $\prob^a(x,D)-\prob^{a'}(x',D)\leq 1-{1\over M}\eta(D)$, so that \eqref{A.1} is satisfied with $\Delta=1-{1\over M}$.
\end{remark}

Using \eqref{A.3} we obtain
\begin{proposition}\label{prop1}
Under \eqref{A.3} we have that
\begin{equation}
\lambda = \sup_{u \ continuous} \lambda^u.
\end{equation}
\end{proposition}
\proof
Following the proof of Lemma 5.1 of \cite{DiMSte1999} one can show that $\lambda$ can be approximated by $\lambda^{u}$, where $u:E\to U$ is piecewise constant and we can assume that points of discontinuity of $u$ are of $\eta$ measure $0$. We now construct a sequence of continuous functions $u_k:E\to U$, $k=1,2,\ldots$ such that $u_k(x)\to u(x)$, as $k\to \infty$ outside of a set of $\eta$ measure $0$. We are going now to show by induction that
whenever bounded sequence of Borel measurable functions $f_k$ converges pointwise to $f$ for each $x\in E$ outside of a set of measure $0$, then for each $m\in \bN$ and for each $x\in E$ outside of a set of measure $0$
\begin{equation}\label{impconv}
\prob_m^{u_k(x)}f_k(x)\to \prob_m^{u(x)}f(x),
\end{equation}
where $\prob_m^{u(x)}f(x)$ stands for $m$-th iteration of the transition operator $\prob^{u(x)}(x,\cdot)$ with $x\in E$.
In fact, for $m=1$ we have
\begin{eqnarray}\label{impconv1}
&&|\prob^{u(x)}f(x)-\prob^{u_k(x)}f_k(x)|\leq |\prob^{u(x)}f(x)-\prob^{u(x)}f_k(x)|+ \nonumber \\
&&|\prob^{u(x)}f_k(x)-\prob^{u_k(x)}f_k(x)|=a_k+b_k
 \end{eqnarray}
and $a_k\to 0$, as $k\to \infty$ by the dominated convergence theorem, while $b_k\to 0$, as $k\to \infty$, outside of a set of $\eta$ measure $0$ by continuity of the transition density $p(x,a,y)$. Consequently we have \eqref{impconv} satisfied for all $x\in E$ outside of a set of $\eta$ measure $0$ with $m=1$. Assume now  $\prob_{m-1}^{u_k(x)}f_k(x)\to \prob_{m-1}^{u(x)}f(x)$, as $k\to \infty$ outside of a set of $\eta$ measure $0$. Then
replacing in \eqref{impconv1} $f$ and $f_k$ by $\prob_{m-1}^{u(x)}f(x)$ and
$\prob_{m-1}^{u_k(x)}f_k(x)$ respectively we obtain the claim for $m$.

Consequently we have
\begin{equation}\label{mainconv}
\int_E c(y,u_k(y))\prob_m^{u_k(x)}(x,dy) \to \int_E c(y,u(y))\prob_m^{u(x)}(x,dy)
\end{equation}
as $k\to \infty$.
Now notice that
\begin{equation}
w^u(x)=\ee_x^u\left[\sum_{i=0}^{n-1}(c(X_i,u(X_i))-\lambda^{u}) + w^{u}(X_n)\right]
\end{equation}
and
\begin{equation}
w^{u_k}(x)=\ee_x^{u_k}\left[\sum_{i=0}^{n-1}(c(X_i,u_k(X_i))-\lambda^{u_k}) + w^{u_k}(X_n)\right].
\end{equation}
Therefore taking into account that we can choose $w^{u_k}$ and $w^u$ such that $\|w^{u_k}\|\leq {1\over 1-\Delta}$ and $\|w^{u}\|\leq {1\over 1-\Delta}$  we have
\begin{equation}
|\lambda^u-\lambda^{u_k}|\leq {1\over n}\mid\ee_x^u\left[\sum_{i=0}^{n-1}(c(X_i,u(X_i))\right]-
\ee_x^{u_k}\left[\sum_{i=0}^{n-1}(c(X_i,u_k(X_i))\right]\mid+{2 \over n(1-\Delta)}
\end{equation}
Letting first $k\to \infty$ (using convergence \eqref{mainconv}) and then $n\to \infty$ we obtain the convergence of $\lambda^{u_k}$ to $\lambda^u$, as $k\to \infty$.
\endproof

\section{Long run risk sensitive generally discounted functional.}\label{Long risk}

Let
\begin{equation}\label{v2}
v^{\gamma,g}(k,x):= \sup_V I_{k,x}^{\gamma,g}(V).
\end{equation}
Bellman equation corresponding to the reward functional \eqref{fun2} depends on the sign of $\gamma$. When $\gamma<0$ we have for $i\geq k$
\begin{equation}\label{B2n}
e^{w^\gamma(i,x)}=\inf_{a\in U}\left[e^{(\phi(i)(c(x,a)-\lambda^\gamma(i)))\gamma}\int_E e^{w^\gamma(i+1,y)}\prob^a(x,dy)\right],
\end{equation}
while for $\gamma>0$
\begin{equation}\label{B2p}
e^{w^\gamma(i,x)}=\sup_{a\in U}\left[e^{(\phi(i)(c(x,a)-\lambda^\gamma(i)))\gamma}\int_E e^{w^\gamma(i+1,y)}\prob^a(x,dy)\right]
\end{equation}
and we are looking for a function $w^\gamma\in C(\bN\times E)$ and a sequence $(\lambda^\gamma(i))$.
Define the following two operators defined for functions $w\in C_{sp}$
\begin{equation}
\Psi^{rg-}w(i,x)=\inf_{a\in U}\ln \left[e^{\phi(i)c(x,a)\gamma}\int_E e^{w(i+1,y)}\prob^a(x,dy)\right]
\end{equation}
and
\begin{equation}
\Psi^{rg+}w(i,x)=\sup_{a\in U}\ln \left[e^{\phi(i)c(x,a)\gamma}\int_E e^{w(i+1,y)}\prob^a(x,dy)\right].
\end{equation}
\begin{proposition}\label{prop2}
Under \eqref{A.1} and \eqref{A.2} the operators $\Psi^{rg-}$ and $\Psi^{rg+}$ are local contractions in the space $C_{sp}$, namely there is a function $\Lambda:(0,\infty)\mapsto (0,1)$ such that whenever $\mid\|w_1\|\mid <M$, $\mid\|w_2\|\mid <M$ we have
\begin{equation}\label{eq21}
\mid \|\Psi^{rgs}w_1-\Psi^{rgs}w_2)\|\mid \leq \Lambda(M)\mid \|w_1-w_2)\|\mid
\end{equation}
with $s\in \left\{-,+\right\}$.
\end{proposition}
\proof
We use the arguments of the paper \cite{DiMSte1999}. For $w\in C_{sp}$ we have \begin{equation}\label{eqpom}
\ln \int_E e^{w(k+1,y)}\prob^a(x,dy)=\sup_{\mu\in {\cal P}(E)}\left[ \int_E w(k+1,y)\mu(dy) - I(\mu, \prob^a(x,dy))\right], \end{equation}
where  $I(\mu, \prob^a(x,dy))$ is the mutual entropy between measures $\mu$ and $\prob^a(x,dy)$.
Moreover we have supremum in \eqref{eqpom} is attained for the measure $\mu\in {\cal P}(E)$ such that for $B\in {\cal E}$
\begin{equation}\label{eqpom1}
\mu(B)= {\int_B e^{w(k+1,z)}\prob^a(x,dz) \over \int_E e^{w(k+1,z)}\prob^a(x,dz)}.
\end{equation}
Using \eqref{eqpom} and \eqref{eqpom1} we obtain that
\begin{eqnarray}\label{eq22}
&&\|\Psi^{rgs}w_1-\Psi^{rgs}w_2\|_{sp}\leq \|w_1(k+1,\cdot)-w_2(k+1,\cdot)\|_{sp} \nonumber \\
&& \sup_{B\in {\cal E}} \sup_{x,x'\in E} \sup_{a,a'\in U} \left(\mu_{k,x,a,w_1}-\mu_{k,x',a',w_2}\right)(B)
\end{eqnarray}
with
\begin{equation}
\mu_{k,x,a,w}(B)={\int_B e^{w(k+1,z)}\prob^a(x,dz) \over \int_E e^{w(k+1,z)}\prob^a(x,dz)}.
\end{equation}
Notice now that
\begin{equation}\label{eqpom2}
e^{-\|w(k+1,\cdot)\|_{sp}}\prob^a(x,B)\leq \mu_{k,x,a,w}(B) \leq e^{\|w(k+1,\cdot)\|_{sp}}\prob^a(x,B).
\end{equation}
We claim that
\begin{eqnarray}\label{eqpom3}
&&\sup_{w_i\in C(E), \|w_i(k+1,\cdot)\|_{sp}\leq M, i=1,2}\  \sup_{B\in {\cal E}} \sup_{x,x'\in E} \sup_{a,a'\in U} \left(\mu_{k,x,a,w_1}-\mu_{k,x',a',w_2}\right)(B):= \nonumber \\
&& \Delta(M)<1.
\end{eqnarray}
In fact, assume contrary to \eqref{eqpom3} that $\mu_{k,x_n,a_n,w_1^n}(B_n)\to 1$ and $\mu_{k,x_n',a_n',w_2^n}(B_n)\to 0$, as $n\to \infty$ for sequences $(w_i^n)$ of functions from $C(\bN\times E)$ such that $\|w_i^n(k+1,\cdot)\|_{sp}\leq M$ for $i=1,2$. Then by \eqref{eqpom2} we have that $\prob^{a_n}(x_n,B_n^c)\to 0$ and $\prob^{a_n'}(x_n',B_n)\to 0$, as $n\to \infty$. Consequently $\prob^{a_n}(x_n,B_n)\to 1$ and $\prob^{a_n'}(x_n',B_n)\to 0$, which contradicts \eqref{A.1}.
From \eqref{eqpom3} together with \eqref{eq22} we obtain \eqref{eq21}, which completes the proof.
\endproof
To get existence of solutions to \eqref{B2n}, \eqref{B2p} we have to know that iterations of the operators $\Phi^{rgs}$ are bounded in the norm $\mid \|\cdot\| \mid$. For this purpose we shall need additional assumptions. We consider first the following
\begin{equation}\label{B.1} \tag{B.1}
\sup_{w\in C(E)} \sup_{x,x'\in E} \sup_{a\in U} {\int_E e^{w(y)}\prob^a(x,dy) \over \int_E e^{w(y)}\prob^{a} (x',dy)}=K<\infty.
\end{equation}
\begin{lemma} \label{lem1} Assumption \eqref{B.1} says that measures $\prob^a(x,dy)$ and $\prob^{a}(x',dy)$ are equivalent with densities uniformly (in $a \in U$, $x,x'\in E$) bounded from below and from above.
\end{lemma}
\proof
Note first that $C(E)$ in \eqref{B.1} may be replaced by $B(E)$ the set of bounded Borel measurable functions. In fact, from \eqref{B.1} it follows that for any $w\in C(E)$
\begin{equation}\label{eq23}
\int_E e^{w(y)}\left(\prob^a(x,dy)-K\prob^{a}(x',dy)\right)\leq 0.
\end{equation}
Denote by $H$ the class of Borel measurable bounded functions for which \eqref{eq23} holds. Clearly $C(E)\subset H$ and bounded limits of functions from $H$ are also in $H$.
Therefore the class $H$ consists of bounded Baire functions, which by Thm 4.5.2 of \cite{Loj1988} coincides with $B(E)$. By Hahn decomposition (see Theorem 7.5.1 of \cite{Loj1988}) there is a set Borel set $D$ such that for each $B\in {\cal E}$ we have that $\prob^a(x,B\cap D)-K \prob^{a}(x',B\cap D)\geq 0$. If $\prob^a(x,D)-K \prob^{a}(x',D)> 0$, then for $w_n(y)=n$ when $y\in D$ and $w_n(y)=0$ otherwise we obtain in \eqref{eq23}
\begin{eqnarray}
&&\int_E e^{w_n(y)}\left(\prob^a(x,dy)-K\prob^{a}(x',dy)\right)= \nonumber \\
&&e^n \left(\prob^a(x,D)-K \prob^{a}(x',D)\right)+1-K-\left(\prob^a(x,D)-K \prob^{a}(x',D)\right)>0
\end{eqnarray}
for a sufficiently large $n$. Consequently we should have $\prob^a(x,D)-K \prob^{a}(x',D)= 0$, and $\prob^a(x,\cdot)-K \prob^{a}(x',\cdot)\leq 0$, which in turn means that $\prob^a(x,\cdot)$ is absolutely continuous with respect to $\prob^{a}(x',\cdot)$ with density dominated by $K$. Since the above consideration holds for any $x,x'\in E$ and $a\in U$ we have the claim of Lemma.
\endproof
Using \eqref{B.1} we obtain
\begin{corollary}\label{cor5}
Under assumption \eqref{B.1} for $w\in C_{sp}$ we have that
\begin{equation}\label{eq23'}
\mid\|\Psi^{rgs}w\|\mid \leq \|\gamma c\|_{sp} + \ln K
\end{equation}
for $s\in \left\{-,+\right\}$.
\end{corollary}
\proof In fact, for $x,x'\in E$ and $k\in \bN$ we have
\begin{eqnarray}
&&\Psi^{rgs}w(k,x)-\Psi^{rgs}w(k,x')\leq \|\gamma c\|_{sp} + \sup_{a\in U}\ln\left({\int_E e^{w(k+1,y)}P^a(x,dy)\over \int_E e^{w(k+1,y)}P^{a}(x',dy)}\right)\leq  \nonumber \\
&&\|\gamma c\|_{sp} + \ln K,
\end{eqnarray}
from which we obtain \eqref{eq23'}.
\endproof
Another condition under which certain iterations of $\Psi^{rgs}$ are bounded is in the form (compare to \cite{DiMSte2000})
\begin{equation}\label{B.2} \tag{B.2}
e^{\|\gamma c\|_{sp}} \Delta <1
\end{equation}
\begin{proposition}\label{prop3}
Under \eqref{A.1}, \eqref{A.2} and \eqref{B.2} for $w_0\equiv 0$ and positive integer $n$ we have the following estimation for the $n$-th iteration of the operator $\Psi^{rgs}$
\begin{equation}\label{eq24}
\|(\Psi^{rgs})^nw_0(k,\cdot)\|_{sp}\leq \|\gamma c \|_{sp} - \ln(1-\Delta e^{\|\gamma c\|_{sp}}):=r(\gamma)
\end{equation}
with $s\in \left\{-,+\right\}$.
\end{proposition}
\proof
For $w\in C_{sp}$ define $\bar{w}(k+1,x)=w(k+1,x)-\inf_{y\in E} w(k+1,y)$. We have
\begin{eqnarray}
&&\Psi^{rgs}w(k,x)-\Psi^{rgs}w(k,x')\leq \|\gamma c\|_{sp} + \sup_{a\in U}\ln\left({\int_E e^{w(k+1,y)}P^a(x,dy)\over \int_E e^{w(k+1,y)}P^{a}(x',dy)}\right)= \nonumber \\
&&\|\gamma c\|_{sp} + \sup_{a\in U}\ln\left({\int_E e^{w(k+1,y)}P^a(x,dy)-\int_E e^{w(k+1,y)}P^{a}(x',dy) \over \int_E e^{w(k+1,y)}P^{a}(x',dy)}+1\right)=\nonumber \\
&&\|\gamma c\|_{sp} + \sup_{a\in U}\ln\left({\int_E e^{\bar{w}(k+1,y)}P^a(x,dy)-\int_E e^{\bar{w}(k+1,y)}P^{a}(x',dy) \over \int_E e^{\bar{w}(k+1,y)}P^{a}(x',dy)}+1\right) \nonumber \\
&&\leq \|\gamma c\|_{sp}+
 \sup_{a\in U}\ln\left( e^{\|{w}(k+1,\cdot)\|_{sp}}(P^a(x,D)-P^{a}(x',D)) + 1\right)\leq \nonumber \\
&&\|\gamma c\|_{sp}+\ln\left( e^{\|{w}(k+1,\cdot)\|_{sp}}\Delta +1\right),
\end{eqnarray}
where $D$ comes from Hahn decomposition of the measure $P^a(x,\cdot)-P^{a}(x',\cdot)$, we used \eqref{A.1} and $\|\bar{w}(k+1,\cdot)\|_{sp}\leq \|{w}(k+1,\cdot)\|_{sp}$.
Consequently
\begin{equation}
\|\Psi^{rgs}w(k,\cdot)\|_{sp} \leq \|\gamma c\|_{sp}+\ln\left( e^{\|{w}(k+1,\cdot)\|_{sp}}\Delta +1\right)
\end{equation}
and
\begin{equation}
\mid \|\Psi^{rgs}w\|\mid \leq \|\gamma c\|_{sp}+\ln\left( e^{\mid\|{w}\|\mid}\Delta +1\right).
\end{equation}
Iterating operator $\Psi^{rgs}$ using the last inequality we obtain
\begin{equation}
\mid \|(\Psi^{rgs})^n w\|\mid \leq \|\gamma c\|_{sp} + \ln\left(\sum_{i=0}^{n-1} (e^{\|\gamma c\|_{sp}}\Delta)^i+ (e^{\|\gamma c\|_{sp}}\Delta)^{n-1}\Delta e^{\mid\|w\|\mid}\right).
\end{equation}
Consequently for $w=w_0\equiv 0$  using  \eqref{B.2} we finally get for each positive integer $n$
\begin{equation}
\mid \|(\Psi^{rgs})^n w_0\|\mid \leq \|\gamma c\|_{sp} - \ln(1-\Delta e^{\|\gamma c\|_{sp}}).
\end{equation}
\endproof

From Proposition \ref{prop2} and Corollary \ref{cor5} or Proposition  \ref{prop3} we have
\begin{theorem}\label{thm3}
Under \eqref{A.1}, \eqref{A.2} and \eqref{B.1} for $\gamma \in \bR$ or under \eqref{A.1}, \eqref{A.2} and \eqref{B.2} for $\gamma$ for which \eqref{B.2} holds we have the existence of function $w^\gamma \in C_{sp}$ such that $\inf_y w^\gamma(i,y)=0$ for $i\geq k$ and sequence $(\lambda^\gamma(i))$, such that equations \eqref{B2n} or \eqref{B2p} are satisfied.
Moreover under \eqref{B.1} we have that $\|w^\gamma\|:=\sup_{x\in E} \sup_{i\geq k} |w^\gamma(i,x)|\leq \|\gamma c\|_{sp}+K$, while under \eqref{B.2} we have that $\|w^\gamma\|\leq r(\gamma)$.   Furthermore
\begin{equation}
v^{\gamma,g}(k,x)=I_{k,x}^{\gamma,g}(\hat{V})=\liminf_{n\to \infty} \left(\sum_{i=k}^{n+k-1}\phi(i)\right)^{-1}\left(\sum_{i=k}^{n+k-1} \phi(i)\lambda^\gamma(i)\right),
\end{equation}
where $\hat{V}=(\hat{u}_k(X_0), \hat{u}_{k+1}(X_1),\ldots, \hat{u}_i(X_i),\ldots)$ is a Markov strategy consisting of the selectors $\hat{u}_i$  to the right hand sides of the equations \eqref{B2n} or \eqref{B2p} respectively.
\end{theorem}
\proof
Under \eqref{B.1}  the operator $\Psi^{rgs}$ is a global contraction in $C_{sp}(E)$. Consequently there is $\bar{w}^\gamma\in C_{sp}(E)$ such that $\mid\|\Psi^{rgs}\bar{w}^\gamma-\bar{w}^\gamma\|\mid=0$. Under \eqref{B.2} the iterations of the operator $\Psi^{rgs}$ starting form $w_0\equiv0$ are bounded and therefore by Proposition \ref{prop2} converge to a fixed point $\bar{w}^\gamma\in C_{sp}(E)$ of the operator $\Psi^{rgs}$. Consequently in both cases there is a sequence $(\bar{\lambda}_i)$ such that $\Psi^{rgs}\bar{w}^\gamma(i,x)-\bar{w}^\gamma(i,x)=\bar{\lambda}_i\gamma$ for $x\in E$ and $i\geq k$. Let ${w}^\gamma(i,x)=\bar{w}^\gamma(i,x)-\inf_{y\in E}\bar{w}^\gamma(i,y)$. Then
\begin{equation}
\Psi^{rgs}{w}^\gamma(i,x)-{w}^\gamma(i,x)=\bar{\lambda}_i\gamma+\inf_{y\in E}\bar{w}^\gamma(i,y)-\inf_{y\in E}\bar{w}^\gamma(i+1,y):=\lambda^\gamma(i)\gamma.
\end{equation}
and consequently the sequence of pairs $(w^\gamma(i,\cdot),\lambda^\gamma(i))$ form a solution to the equations \eqref{B2n} or \eqref{B2p} respectively. The estimates for $w^\gamma$ follow from the proofs of Corollary \ref{cor5} and Proposition \ref{prop3}.
\endproof
In the case of long run risk sensitive functional we have the following two Bellman equations depending on the sign of $\gamma$. We are looking for a function $w\in C(E)$ and a constant $\lambda$ such that when $\gamma<0$
\begin{equation}\label{B4n}
e^{w^\gamma(x)}=\inf_{a\in U}\left[e^{(c(x,a)-\lambda^\gamma)\gamma}\int_E e^{w^\gamma(y)}\prob^a(x,dy)\right],
\end{equation}
while when $\gamma>0$
\begin{equation}\label{B4p}
e^{w^\gamma(x)}=\sup_{a\in U}\left[e^{(c(x,a)-\lambda^\gamma)\gamma}\int_E e^{w^\gamma(y)}\prob^a(x,dy)\right].
\end{equation}

Following proofs of Proposition \ref{prop2}, Corollary \ref{cor5}, Proposition  \ref{prop3} and Theorem \ref{thm3} we have

\begin{corollary}\label{cor6}
Under \eqref{A.1}, \eqref{A.2} and \eqref{B.1} for $\gamma\in \bR$ or under \eqref{A.1}, \eqref{A.2} and \eqref{B.2} for $\gamma$ for which \eqref{B.2} holds there is a function $w^\gamma\in C(E)$ such that $\inf_{y\in E}w^\gamma(y)=0$ and a constant $\lambda^\gamma$ which are solutions to the equations \eqref{B4n} or \eqref{B4p}. Furthermore under \eqref{B.1} we have that $\|w^\gamma\|:=\sup_{x\in E} \leq \|\gamma c\|_{sp}+K$, while under \eqref{B.2} we have that $\|w^\gamma\|\leq r(\gamma)$.
Optimal value of the reward functional \eqref{fun4} is $\lambda^\gamma$ and optimal strategies are stationary $a_i=\hat{u}^\gamma(X_i)$, where $\hat{u}^\gamma$ is a selector of the right hand side of the equation \eqref{B4n} for negative $\gamma$ or \eqref{B4p} for positive $\gamma$.
\end{corollary}

By analogy to Section \ref{Average} we consider now the problems with fixed Markov controls $V_m^k=(u_k(X_0),u_{k+1}(X_1),\ldots,u_{k+i}(X_i),\ldots)$.

\begin{proposition}\label{prop4} Under \eqref{A.1}, \eqref{A.2} and \eqref{B.1} for $\gamma\in \bR$ or under \eqref{A.1}, \eqref{A.2} and \eqref{B.2} for $\gamma$ for which \eqref{B.2} holds and Markov control $V_m^k$ there is a function $w^{V_m^k,\gamma}\in B_{sp}$ such that $\inf_{y\in E}w^{V_m^k,\gamma}(i,y)=0$ for $i\geq k$  and sequence $(\lambda^{V_m^k,\gamma}(i))$ such that the following multiplicative Poisson equation is satisfied for $i\geq k$
\begin{equation}\label{P2}
e^{w^{V_m^k,\gamma}(i,x)}=e^{(\phi(i)(c(x,u_i(x))-\lambda^{V_m^k,\gamma}(i)))\gamma}\int_E e^{w^{V_m^k,\gamma}(i+1,y)}\prob^{u_i(x)}(x,dy).
\end{equation}
Moereover under \eqref{B.1} we have that $\|w^{V_m^k,\gamma}\|:=\sup_{x\in E} \sup_{i\geq k} |w^{V_m^k,\gamma}(i,x)|\leq \|\gamma c\|_{sp}+K$, while under \eqref{B.2} we have that $\|w^{V_m^k,\gamma}\|\leq r(\gamma)$.
\end{proposition}
\proof We follow first the proof of Proposition \ref{prop2} for the operator
\begin{equation}
\Psi^\gamma w(i,x)=\ln \left[e^{\phi(i)c(x,u_i(x))\gamma}\int_E e^{w(i+1,y)}\prob^{u_i(x)}(x,dy)\right]
\end{equation}
defined for functions $w\in B_{sp}$ with $i\geq k$. By similar consideration as in the proof of Proposition \ref{prop2} it is a local contraction. By Corollary \ref{cor5} under \eqref{B.1}
we have
\begin{equation}
\mid\|\Psi_{\gamma} w\|\mid \leq \|\gamma c\|_{sp} + \ln K,
\end{equation}
while under \eqref{B.2}  for $w_0\equiv 0$ and positive integer $n$, by Proposition \ref{prop3} (for $\gamma$ such that \eqref{B.2} holds)
\begin{equation}
\|(\Psi^\gamma)^n w_0(k,\cdot)\|_{sp}\leq \|\gamma c \|_{sp} - \ln(1-\Delta e^{\|\gamma c\|_{sp}})=r(\gamma)
\end{equation}
The remaining part of the proof follows from Theorem \ref{thm3}.
\endproof

In the case of Markov stationary control $V_m=(u(X_0),u(X_1),\ldots, u(X_i),\ldots):=u$, where $u:E\to U$ is Borel measurable and absence of discount factor $(\phi(i))$  as in Proposition \ref{prop4} we obtain
\begin{corollary}\label{cor7}
Under \eqref{A.1}, \eqref{A.2} and \eqref{B.1} for $\gamma\in \bR$ or under \eqref{A.1}, \eqref{A.2} and \eqref{B.2} for $\gamma$ for which \eqref{B.2} holds,
  for stationary Markov control $V_m=u$  there is a bounded function $w^{u,\gamma}\in B(E)$ such that $\inf_{y\in E} w^{u,\gamma}(y)= 0$ and constant $\lambda^{u,\gamma}$ such that the multiplicative Poisson equation
\begin{equation}\label{P4}
e^{w^{u,\gamma}(x)}=e^{(c(x,u(x))-\lambda^{u,\gamma})\gamma}\int_E e^{w^{u,\gamma}(y)}\prob^{u(x)}(x,dy)
\end{equation}
is satisfied. Moreover under \eqref{B.1} we have that $\|w^{u,\gamma}\|:=\sup_{x\in E} |w^{u,\gamma}(x)|\leq \|\gamma c\|_{sp}+K$, while under \eqref{B.2} we have that $\|w^{u,\gamma}\|\leq r(\gamma)$ and
\begin{equation}\label{vvvp4}
\lambda^{u,\gamma}=I_x^\gamma(u)=\lim_{n\to \infty}{1\over n \gamma} \ln \left(\ee_x^u\left[\exp\left\{\gamma\sum_{i=0}^{n-1} c(X_i,u(X_i))\right\}\right]\right).
\end{equation}
\end{corollary}
By analogy to Proposition \ref{prop1} we have the following result
\begin{proposition}\label{prop5}
Under \eqref{A.3}  for $\gamma\in \bR$  we have
\begin{equation}
\lambda^{\gamma}=\sup_{u \  continuous} \lambda^{u,\gamma}.
\end{equation}
\end{proposition}
\proof
We adapt the proof of Lemma 5.1 of \cite{DiMSte1999} to approximate $\lambda^{\gamma}$ by $\lambda^{u,\gamma}$ with piecewise constant $u:E\to U$ with discontinuity points of measure $\eta$ equal to $0$. As in Proposition \ref{prop1} we have that for any
bounded sequence of Borel measurable functions $f_k$ converging pointwise to $f$ for each $x\in E$ outside of the set of $\eta$ measure $0$,   we have
$\prob^{u_k(x)}f_k(x)\to \prob^{u(x)}f(x)$ as $k\to \infty$ for each $x$ outside of the set of $\eta$ measure $0$.
Therefore
\begin{equation}
e^{w^{u,\gamma}(x)}=\ee_x^u\left[\exp\left(\sum_{i=0}^{n-1}\gamma(c(X_i,u(X_i))-\lambda^{u,\gamma}) + w^{u,\gamma}(X_n)\right)\right]
\end{equation}
and
\begin{equation}
e^{w^{u_k,\gamma}(x)}=\ee_x^{u_k}\left[\exp\left(\sum_{i=0}^{n-1}\gamma(c(X_i,u_k(X_i))-\lambda^{u_k,\gamma}) + w^{u_k,\gamma}(X_n)\right)\right].
\end{equation}
Finally
\begin{eqnarray}
&&|\gamma(\lambda^{u,\gamma}-\lambda^{u_k,\gamma})|\leq {1\over n}\left[ \ln\left({\ee_x^u\left[\exp\left(\sum_{i=0}^{n-1}\gamma c(X_i,u(X_i)\right)\right]\over
\ee_x^{u_k}\left[\exp\left(\sum_{i=0}^{n-1}\gamma c(X_i,u_k(X_i)\right)\right]}\right)+ \right. \nonumber \\
&&\left. \|w^{u,\gamma}\|+\|w^{u_k,\gamma}\| \right]
\end{eqnarray}
and since the norms of $w^{u_k,\gamma}$ are bounded letting first $k\to \infty$ and then $n\to \infty$ we obtain that $\lambda^{u_k,\gamma}\to \lambda^{u,\gamma}$ as $k\to \infty$.
\endproof

By analogy to Section \ref{Average} we want to compare optimal values of the functional \eqref{fun2} and \eqref{fun4}.

\begin{theorem}\label{thm4}
Under \eqref{A.1}, \eqref{A.2} and \eqref{B.1} for $\gamma>0$ or \eqref{A.1}, \eqref{A.2} and  \eqref{B.2} for $\gamma>0$ for which \eqref{B.2} holds we have
\begin{equation}\label{posrisk}
\sup_V I_{k,x}^{\gamma,g}(V) = v^{\gamma,g}(k,x)\leq\lambda^\gamma
\end{equation}
\end{theorem}
\proof
From Bellman equation \eqref{B4p} for any strategy $V=(a_0,a_1,\ldots,a_n,\ldots)$ we have
\begin{equation}
e^{\gamma c(X_i,a_i)}\leq e^{w^\gamma(X_i)+\lambda^\gamma \gamma}{1\over \ee^V\left[e^{w^\gamma(X_{i+1})}|F_i\right]}
\end{equation}
and therefore
\begin{eqnarray}\label{innq0}
&&\exp\left(\gamma \sum_{i=k}^{n+k-1}\phi(i)c(X_{i-k},a_{i-k})\right)\leq \nonumber \\
&&e^{\lambda^\gamma \gamma \sum_{i=k}^{n+k-1}\phi(i)} \prod_{j=k}^{n+k-1} {e^{\phi(j)w^\gamma(X_{j-k})}\over \ee_{x}^V\left[e^{w^\gamma(X_{j+1-k})}|F_{j-k}\right]^{\phi(j)}}= \nonumber \\
&&e^{\lambda^\gamma \gamma \sum_{i=k}^{n+k-1}\phi(i)}e^{\phi(k)w^\gamma(x)} {1 \over
\ee_{x}^V\left[e^{w^\gamma(X_{n})}|F_{n-1}\right]^{\phi(n+k-1)}} \nonumber \\
&& \prod_{j=k}^{n+k-2} {e^{\phi(j+1)w^\gamma(X_{j+1-k})}\over \ee_{x}^V\left[e^{w^\gamma(X_{j+1-k})}|F_{j-k}\right]^{\phi(j)}}.
\end{eqnarray}
Consequently
\begin{eqnarray} \label{innq0'}
&&\ee_x^V\left[\exp\left(\gamma \sum_{i=k}^{n+k-1}\phi(i)c(X_{i-k},a_{i-k})\right)\right]\leq \nonumber \\
&&e^{\lambda^\gamma \gamma \sum_{i=k}^{n+k-1}\phi(i)}e^{\phi(k)w^\gamma(x)}
\ee_x^V\left[ {1 \over
\ee_{x}^V\left[e^{w^\gamma(X_{n})}|F_{n-1}\right]^{\phi(n+k-1)}}\right. \nonumber \\
&&\left. \prod_{j=k}^{n+k-2} {e^{\phi(j+1)w^\gamma(X_{j-k+1})}\over \ee_{x}^V\left[e^{w^\gamma(X_{j-k+1})}|F_{j-k}\right]^{\phi(j)}}\right]
\end{eqnarray}
and assuming that $w^\gamma\geq 0$ (otherwise we consider solution  to \eqref{B4p} of the form   $w^\gamma-\inf_{y\in E}w^\gamma(y)$ with the same $\lambda^\gamma$)
\begin{eqnarray}\label{innq1}
&&\ee_x^V\left[ {1 \over
\ee_{x}^V\left[e^{w^\gamma(X_{n})}|F_{n-1}\right]^{\phi(n+k-1)}}\prod_{j=k}^{n+k-2} {e^{\phi(j+1)w^\gamma(X_{j-k+1})}\over \ee_{x}^V\left[e^{w^\gamma(X_{j-k+1})}|F_{j-k}\right]^{\phi(j)}}\right]\leq \nonumber \\
&&\ee_x^V\left[
\prod_{j=k}^{n+k-3} {e^{\phi(j+1)w^\gamma(X_{j-k+1})}\over \ee_{x}^V\left[e^{w^\gamma(X_{j-k+1})}|F_{j-k}\right]^{\phi(j)}}
{\ee_x^V\left[e^{\phi(n+k-1)w^\gamma(X_{n-1})}|F_{n-2}\right]\over \ee_{x}^V\left[e^{w^\gamma(X_{n-1})}|F_{n-2}\right]^{\phi(n+k-2)}}\right].
\end{eqnarray}
Since by conditional H\"older inequality
\begin{eqnarray}\label{innq2}
&&{\ee_x^V\left[e^{\phi(n+k-1)w^\gamma(X_{n-1})}|F_{n-2}\right]\over \ee_{x}^V\left[e^{w^\gamma(X_{n-1})}|F_{n-2}\right]^{\phi(n+k-2)}}\leq \nonumber \\ &&\ee_{x}^V\left[e^{w^\gamma(X_{n-1})}|F_{n-2}\right]^{\phi(n+k-1)-\phi(n+k-2)}\leq e^{-\|w^\gamma\|(\phi(n+k-1)-\phi(n+k-2))}
\end{eqnarray}
using succesively \eqref{innq2} in  \eqref{innq1} and then substituting it to \eqref{innq0'} we obtain
\begin{eqnarray}\label{innq3}
&&\ee_x^V\left[\exp\left(\gamma \sum_{i=k}^{n+k-1}\phi(i)c(X_{i-k},a_{i-k})\right)\right]\leq \nonumber \\
&&e^{\lambda^\gamma \gamma \sum_{i=k}^{n+k-1}\phi(i)}e^{\phi(k)w^\gamma(x)}
e^{\|w^\gamma\|(\phi(k)-\phi(n+k-1))}
\end{eqnarray}
Finally taking logarithm in both sides of \eqref{innq3}, then dividing them by $\gamma\sum_{i=k}^{n+k-1}\phi(i)$ and letting $n\to \infty$ we obtain that $I_{k,x}^{\gamma,g}(V)\leq \lambda^\gamma$.
\endproof
Following the proof of the last theorem we easily obtain
\begin{corollary}\label{cor8}
Under assumptions of Theorem \ref{thm4} for $\gamma>0$ and stationary Markov control $u:E\to U$ we have
\begin{equation}
I_{k,x}^{\gamma,g}(u)\leq \lambda^{u,\gamma}
\end{equation}
\end{corollary}
\begin{remark}\label{rem5}
For an optimal Markov strategy $\hat{V}=\hat{u}$ for the functional \eqref{fun4} from \eqref{innq0} we have
\begin{equation}
\ee_x^{\hat{V}}\left[\exp\left(\gamma \sum_{i=k}^{n+k-1}\phi(i)c(X_{i-k},\hat{u}(X_{i-k})\right)\right]=
e^{\lambda^\gamma \gamma \sum_{i=k}^{n+k-1}\phi(i)}e^{\phi(k)w^\gamma(x)}d_k^n,
\end{equation}
where
\begin{equation}
d_k^n:=\ee_x^{\hat{V}}\left[ {1 \over
\ee_{x}^{\hat{V}}\left[e^{w^\gamma(X_{n})}|F_{n-1}\right]^{\phi(n+k-1)}}\prod_{j=k}^{n+k-2} {e^{\phi(j+1)w^\gamma(X_{j-k+1})}\over \ee_{x}^{\hat{V}}\left[e^{w^\gamma(X_{j-k+1})}|F_{j-k}\right]^{\phi(j)}}\right].
\end{equation}
To get equality in \eqref{posrisk} we need an estimation from below of
\begin{equation}
{\ee_x^{\hat{V}}\left[e^{\phi(n+k-1)w^\gamma(X_{n-1})}|F_{n-2}\right]\over \ee_{x}^{\hat{V}}\left[e^{w^\gamma(X_{n-1})}|F_{n-2}\right]^{\phi(n+k-2)}}.
\end{equation}
which does not seem to be easy.
\end{remark}

\section{Risk sensitive asymptotics.}\label{Risk asympt}

For any bounded random variable $Y$, by H\"older and Jensen inequalities, for $\gamma>0$ we have
\begin{equation}
{1\over -\gamma}\ln\left( \ee\left[e^{-\gamma Y}\right]\right)\leq \ee[Y] \leq {1\over \gamma}\ln\left(\ee\left[e^{\gamma Y}\right]\right).
\end{equation}
Therefore for any control $V$ and $\gamma>0$ we have
\begin{equation}\label{iiwi}
I_{k,x}^{-\gamma,g}(V)\leq J_{k,x}^g(V) \leq I_{k,x}^{\gamma,g}(V).
\end{equation}
Consequently by Theorem \ref{thm2} for $\gamma >0$
\begin{equation}\label{crineq}
v^{-\gamma,g}(k,x)\leq \lambda \leq v^{\gamma,g}(k,x).
\end{equation}
We have even more
\begin{theorem}\label{thm5}
Under \eqref{A.3} we have that $\lim_{\gamma\to 0^-}v^{\gamma,g}(k,x)=\lambda$, for any $k\in \bN$ and $x\in E$. Furthermore almost optimal continuous stationary Markov control ${u}$ for the functional \eqref{fun3} is also nearly optimal for functional \eqref{fun2} provided that $\gamma<0$ is sufficiently close to $0$. In the case when $ \gamma>0$ we have that for each continuous stationary Markov control $u:E \to U$
\begin{equation}\label{pconver}
\lim_{\gamma\to 0} I_{k,x}^{\gamma,g}(u)=\lim_{\gamma\to 0}I_x^\gamma(u)=\lambda^u.
\end{equation}
\end{theorem}
\proof
By Proposition \ref{prop1} for a given $\epsilon>0$ there is a continuous $u:E\to U$ such that $\lambda^u\geq \lambda-\epsilon$. Markov process $(X_n)$ controlled with $V=\left\{u(X_0),u(X_1),\ldots, u(X_i),\ldots\right\}$ is Feller for which we can use large deviation estimates formulated in the Appendix. By \eqref{A.3} it has a unique invariant measure $\mu$ and for the transition operator $\prob^{u(x)}(x,\cdot)$ assumption \eqref{D.1} is satisfied. Therefore for $\epsilon>0$ the set $C=\left\{\nu\in {\cal P}(E): |\int_E c(y,u(y))\nu(dy)-\int_E c(y,u(y))\mu(dy)|\geq \epsilon\right\}$ is closed and by Lemma \ref{alem3} $\inf_{\zeta \in C} I(\zeta)>0$.
Consequently by Theorem \ref{athm1} there is $e>0$ and $N$ such that for $n\geq N$ we have $\prob_{k,x}\left(A_\epsilon^n\right)\leq \exp(-\sum_{i=k}^{n+m-1} \phi(i)e)$, where
$A_\epsilon^n=\left\{\omega: |\left(\sum_{i=k}^{n+m-1} \phi(i)\right)^{-1}\sum_{i=k}^{n+k-1} \phi(i)c(X_{i-k},u(X_{i-k}))-\int_E c(y,u(y))\mu(dy)|\geq \epsilon\right\}$.
Therefore for any $\gamma\neq 0$ such that $|\gamma|<{e\over 2\|c\|}$  we have
\begin{eqnarray}\label{iiin}
&& \ln\left( \ee_{k,x}^V\left[\exp\left(\gamma\sum_{i=k}^{n+k-1}\phi(i)c(X_{i-k},u(X_{i-k}))
\right)\right]\right) = \nonumber \\
&& \int_E c(y,u(y))\mu(dy)\gamma\sum_{i=k}^{n+k-1}\phi(i)+\nonumber \\
 &&\ln\left( \ee_{k,x}^V\left[\exp\left(\gamma\sum_{i=k}^{n+k-1}\phi(i)\left(c(X_{i-k},u(X_{i-k}))-\int_E c(y,u(y))\mu(dy)\right)\right)\right]\right)\leq \nonumber \\
&&\int_E c(y,u(y))\mu(dy)\gamma\sum_{i=k}^{n+k-1}\phi(i)+ \nonumber \\
&& \ln\left( \ee_{k,x}^V\left[1_{A_\epsilon^n}\exp\left(|\gamma|\sum_{i=k}^{n+k-1}\phi(i)2\|c\|\right)
+1_{\Omega\setminus A_\epsilon^n}\exp\left(|\gamma|\sum_{i=k}^{n+k-1}\phi(i)\epsilon
\right) \right]\right)\leq \nonumber \\
&&\int_E c(y,u(y))\mu(dy)\gamma\sum_{i=k}^{n+k-1}\phi(i)+ \nonumber \\
&&\ln\left(
\exp\left(|\gamma|\sum_{i=k}^{n+k-1}\phi(i)(2\|c\|-{e\over |\gamma|})\right)+
\exp\left(|\gamma|\sum_{i=k}^{n+k-1}\phi(i)\epsilon\right)\right)\leq \nonumber \\
&&\int_E c(y,u(y))\mu(dy)\gamma\sum_{i=k}^{n+k-1}\phi(i)+ \ln\left(1+
\exp\left(|\gamma|\sum_{i=k}^{n+k-1}\phi(i)\epsilon\right)\right)\leq \nonumber \\
&&\int_E c(y,u(y))\mu(dy)\gamma\sum_{i=k}^{n+k-1}\phi(i)+|\gamma|\sum_{i=k}^{n+k-1}\phi(i)\epsilon.
\end{eqnarray}
Consequently for $\gamma<0$ such that $|\gamma|<{e\over 2\|c\|}$  we have that
\begin{equation}
I_{k,x}^{\gamma,g}(V)\geq \int_E c(y,u(y))\mu(dy)- \epsilon=\lambda^u-\epsilon\geq \lambda-2\epsilon
\end{equation}
and by \eqref{crineq} $I_{k,x}^{\gamma,g}(V)\geq v^{\gamma,g}(k,x)-2\epsilon$, which means that control $V$ is $2\epsilon$ - optimal for $I_{k,x}^{\gamma,g}$. Since $I_{k,x}^{\gamma,g}(V)\leq v^{\gamma,g}(k,x)\leq \lambda$ and $\epsilon$ could be chosen arbitrarily small we have also convergence of $v^{\gamma,g}(k,x)$ to $\lambda$, as $\gamma\to 0$.

In the case  $\gamma>0$ consider continuous Markov control $u:E\to E$. Then for $\phi(i)\equiv 1$ following \eqref{iiin} for strategy $V=\left\{u(X_0),u(X_1),\ldots, u(X_i),\ldots\right\}$ we obtain for $0<\gamma<{e\over 2\|c\|}$
\begin{equation}
I_x^\gamma(V)=\lambda^{u,\gamma}\leq \lambda^u+\epsilon.
\end{equation}
Since by Corollary \ref{cor8} and \eqref{iiwi} we have that $\lambda^{u,\gamma}\geq I_{k,x}^{\gamma,g}(u)\geq \lambda^u$ we obtain \eqref{pconver}, which completes the proof.

\endproof

\begin{remark}
In the case of positive risk factor $\gamma$ we have definitely weaker result. Namely we have convergence \eqref{pconver} for each fixed continuous Markov control. We could get a similar result to that for negative $\gamma$ if we were able to show that
\begin{equation}\label{row0}
\sup_u \lambda^{u,\gamma}\to \sup_u \lambda^u,
\end{equation}
as $\gamma \to 0^+$, with supremum over continuous Markov controls $u$.
By Theorem 6.9 and Corollary 7.21 of \cite{Stroock} under \eqref{A.3} we obtain that
\begin{equation}\label{row1}
I_x^\gamma(u)\leq \sup_{\nu\in {\cal P}(E)} \left[\int_E c(z,u(z))\nu(dz)
-{1\over \gamma} I^u(\nu)\right],
\end{equation}
where
\begin{equation}\label{row2}
I^u(\nu)=\sup_{f\in F} \int_E \ln {f(y) \over \prob^u f(y)} \nu(dy)
\end{equation}
and $F$ is the class is continuous positive bounded functions $f$ which are separated from $0$.
For a given $\epsilon>0$  there are a continuous Markov control $u_\gamma$ and measure $\nu^{u_\gamma,\gamma}$ such that
\begin{equation}\label{row3}
\sup_u\sup_{\nu\in {\cal P}(E)} \left[\int_E c(z,u(z))\nu(dz)
-{1\over \gamma} I^u(\nu)\right]\leq \epsilon + \int_E c(z,u_\gamma(z))\nu^{u_\gamma,\gamma}(dz)
-{1\over \gamma} I^{u_\gamma}(\nu^{u_\gamma,\gamma}).
\end{equation}
Letting $\gamma\to 0^+$ we see that $I^{u_\gamma}(\nu^{u_\gamma,\gamma})\to 0$. When $E$ is compact we can find a subsequence $\gamma_k \to 0^+$ and measure ${\bar{\nu}}\in {\cal P}(E)$ such that $\nu^{u_{\gamma_k},\gamma_k} \rightarrow {\bar{\nu}}$, as $k\to \infty$. The problem is  that only in the case of denumerable state space we can expect that certain subsequence, for simplicity still denoted by $u_{\gamma_k}$, converges to a function $u$ uniformly (or uniformly on compact sets) and then taking into account lower semicontinuity of $\nu \to I^u(\nu)$ we obtain that $I^u(\bar{\nu})=0$. Therefore $\bar{\nu}$ is an invariant measure for $\prob^{u(x)}(x,\cdot)$. Consequently from \eqref{row1} and \eqref{row3} we obtain
\begin{equation}
\limsup_{\gamma \to 0} \sup_u I_x^\gamma(u)-\epsilon\leq \int_E c(z,u(z)){\bar\nu(dz)}=\lambda^u,
\end{equation}
which, since $\epsilon$ could be chosen arbitrarily small, completes the proof of \eqref{row0} when $E$ is finite.
\end{remark}


\section{Appendix}\label{secA1}

We formulate first and prove a large deviation result concerning upper estimates for discounted empirical measures of Feller Markov processes $(X_n)$ with transition operator $\prob(x,\cdot)$ such that $\prob C(E)\subset C(E)$.  Denote by ${\cal P}(E)$ and ${\cal P}({\cal P}(E))$ the sets of probability measures on $E$ or on ${\cal P}(E)$ the set of probability measures on $E$ endowed with weak convergence topology.
Let for $B\in {\cal E}$, $k\in \bN$, $x\in E$
\begin{equation}
S_{k}^n(B):= \left(\sum_{i=k}^{n+k-1} \phi(i)\right)^{-1} \sum_{i=k}^{n+k-1} \phi(i) 1_B(X_{i-k}).
\end{equation}
Clearly $S_{k}^n\in {\cal P}(E)$. For a Borel measurable subset $D\subset {\cal P}(E)$ define
\begin{equation}
Q_{k,x}^n(D):= \prob_{k,x}\left\{S_{k}^n\in D\right\}.
\end{equation}
Define the sets of functions $F:=\left\{f\in C(E): \exists_{a>0} f(x)\geq a, \ \  for \ \ x\in E\right\}$ and $F_d:=\left\{f\in F: \sup_{x\in E} f(x) \leq d \inf_{x\in E} f(x)\right\}$ for $d>0$.
For $\nu\in {\cal P}(E)$ let
\begin{equation}
I(\nu):=\sup_{f\in F} \int_E \ln {f(x) \over \prob f(x)} \nu(dx)
\end{equation}
and for $d>0$
\begin{equation}
I_d(\nu):=\sup_{f\in F_d} \int_E \ln {f(x) \over \prob f(x)} \nu(dx).
\end{equation}
One can notice that because of the form of $I$ and $I_d$ we can restrict ourselves to functions $f\in F$ such that $f(x)\geq 1$ for $x\in E$.
Consider the following assumption:
there is $\eta\in {\cal P}(E)$ such that
\begin{equation}\label{D.1} \tag{D.1}
\prob(x,B)=\int_B p(x,y)\eta(dy) \ \ for \ \ B\in{\cal E} \ and  \ x\in E,
\end{equation}
the mapping $E\ni x \to p(x,y)$ is continuous for $y\in E$ and there is $M>0$ such that ${1\over M} \leq p(x,y) \leq  M$ for $x,y\in E$.
Then we have
\begin{lemma}\label{alem1}
When the set $E$ is compact or $E$ is complete separable and \eqref{D.1} is satisfied then  the set
$C_q=\left\{\zeta\in {\cal P}(E): I(\zeta)\leq q\right\}$ is compact in ${\cal P}(E)$ for each $q>0$.
\end{lemma}
\proof Following the proof Theorem 5.2 of \cite{DiMSte1999} we show that condition $H^*$ of \cite{DV3} is satisfied and therefore by Lemma 4.2 of \cite{DV3} we obtain compactness of $C_q$.
\endproof

We have the following upper estimate, which generalizes section 2 of \cite{DV1}
\begin{theorem}\label{athm1}
For a compact set $C \subset {\cal P}(E) $ or closed set $C \subset {\cal P}(E) $ under additional assumption \eqref{D.1} we have
\begin{equation}
\limsup_{n\to \infty} \left(\sum_{i=k}^{n+k-1} \phi(i)\right)^{-1} \sup_{x\in E} \ln Q_{k,x}^n(C)\leq - \inf_{\nu \in C} I(\nu).
\end{equation}
\end{theorem}
\proof For $f\in F_d$ such that $\inf_{y\in E} f(y)\geq 1$, taking into account that $\phi(i)\leq 1$, we have
\begin{eqnarray}\label{aa1}
&& \ee_{k,x}\left[ \exp \left\{ \sum_{i=k}^{n+k-1} \phi(i) \ln {f(X_{i-k})\over \prob f(X_{i-k})}\right\}\right]= \nonumber \\
&& \ee_{k,x}\left[{f(x)^{\phi(k)}\over (\prob f(X_{n-1}))^{\phi(n+k-1)}} \prod_{i=k}^{n+k-2} {f(X_{i-k+1})^{\phi(i+1)} \over (\prob f(X_{i-k}))^{\phi(i)}}\right]\leq  \nonumber \\
&& d \ee_{k,x}\left[\prod_{i=k}^{n+k-2} {f(X_{i-k+1})^{\phi(i+1)} \over (\prob f(X_{i-k}))^{\phi(i)}}\right]
\end{eqnarray}
Now, using H\"older inequality since $\phi(n+k-1) \leq \phi(n+k-2)$ we obtain
\begin{eqnarray}\label{aa2}
&&\ee_{k,x}\left[\prod_{i=k}^{n+k-2} {f(X_{i-k+1})^{\phi(i+1)} \over (\prob f(X_{i-k}))^{\phi(i)}}\right]=  \nonumber \\
&& \ee_{k,x}\left[\prod_{i=k}^{n+k-3} {f(X_{i-k+1})^{\phi(i+1)} \over (\prob f(X_{i-k}))^{\phi(i)}}\ee_{k,x}\left[{f(X_{n-1})^{\phi(n+k-1)} \over (\prob f(X_{n-2}))^{\phi(n+k-2)}}|F_{n-2}\right]\right]= \nonumber \\
&& \ee_{k,x}\left[\prod_{i=k}^{n+k-3} {f(X_{i-k+1})^{\phi(i+1)} \over (\prob f(X_{i-k}))^{\phi(i)}}{(\prob f(X_{n-2}))^{\phi(n+k-1)}\over (\prob f(X_{n-2}))^{\phi(n+k-2)}}\right]\leq \nonumber \\
&& \ee_{k,x}\left[\prod_{i=k}^{n+k-3} {f(X_{i})^{\phi(i+1)} \over (\prob f(X_{i-k}))^{\phi(i)}}\right].
\end{eqnarray}
Therefore iterating \eqref{aa2} we get
\begin{equation}\label{aa3}
\ee_{k,x}\left[\prod_{i=k}^{n+k-2} {f(X_{i-k+1})^{\phi(i+1)} \over (\prob f(X_{i-k}))^{\phi(i)}}\right]\leq 1
\end{equation}
and consequently using \eqref{aa1} we obtain
\begin{equation}\label{aa4}
\ee_{k,x}\left[ \exp \left\{ \sum_{i=k}^{n+k-1} \phi(i) \ln {f(X_{i-k})\over \prob f(X_{i-k})}\right\}\right]\leq d.
\end{equation}
We can rewrite \eqref{aa4} in the form
\begin{equation}\label{aa5}
\int_{{\cal P}(E)} \exp\left\{\sum_{i=k}^{n+k-1} \phi(i) \int_E \ln {f(y)\over \prob f(y)}\zeta(dy)\right\} Q_{k,x}^n(d\zeta)\leq d.
\end{equation}
For Borel measurable subset $B\subset {\cal P}(E)$
we have
\begin{equation}\label{aa6}
Q_{k,x}^n(B)\leq d \exp\left\{-\sum_{i=k}^{n+k-1} \phi(i) \inf_{\zeta\in B} \int_E \ln {f(y)\over \prob f(y)} \zeta(dy)\right\}.
\end{equation}
Let $C\subset {\cal P}(E)$ be a compact set and
\begin{equation}\label{aa7}
\inf_{\zeta \in C} I_d(\zeta)>\kappa>0.
\end{equation}
By the definition of $I_d$ we have that $C\subset \cup_{f\in F_d} \left\{\zeta \in {\cal P}(E): \zeta(\ln {f\over \prob f})>\kappa \right\}$. Since $C$ is compact there is a finite set $\left\{f_1,f_2,\ldots,f_m\right\}\subset F_d$ such that $C\subset \cup_{j=1}^m \left\{\zeta \in {\cal P}(E): \zeta(\ln {f_j\over \prob f_j})>\kappa \right\}$. For the sets $K_j:=\left\{\zeta \in {\cal P}(E): \zeta(\ln {f_j\over \prob f_j})>\kappa \right\}\cap C$ from \eqref{aa6} and \eqref{aa7} we obtain
\begin{equation}\label{aa8}
Q_{k,x}^n(K_j)\leq d \exp\left\{-\sum_{i=k}^{n+k-1} \phi(i) \inf_{\zeta\in K_j} \int_E \ln {f(y)\over \prob f(y)} \zeta(dy)\right\}\leq d \exp\left\{-\sum_{i=k}^{n+k-1}\phi(i)\kappa\right\}.
\end{equation}
Therefore
\begin{equation}\label{aa9}
Q_{k,x}^n(C) \leq \sum_{j=1}^m Q_{k,x}^n(K_j)\leq m d \exp\left\{-\sum_{i=k}^{n+k-1}\phi(i)\kappa\right\}.
\end{equation}
Since $\kappa$ could be chosen arbitrarily close to $\inf_{\zeta \in C} I_d(\zeta)$ we have
\begin{equation}\label{aa10}
\limsup_{n\to \infty} \sup_{x\in E}(\sum_{i=k}^{n+k-1}\phi(i))^{-1}\ln Q_{k,x}^n(C)\leq - \inf_{\zeta \in C} I_d(\zeta).
\end{equation}
Similarly as above we have that $C\subset \cup_{f\in F} \left\{\zeta \in {\cal P}(E): \zeta(\ln {f\over \prob f})>  \inf_{\zeta \in C} I_d(\zeta) - \epsilon \right\}$ for $\epsilon>0$ and therefore by compactness of $C$ there is a finite set $\left\{f_1,\ldots,f_p\right\}\subset F$ such that $C\subset \cup_{i=1}^p \left\{\zeta \in {\cal P}(E): \zeta(\ln {f_i\over \prob f_i})>  \inf_{\zeta \in C} I_d(\zeta) - \epsilon \right\}$. Therefore there is $d>0$ such that $\left\{f_1,\ldots,f_p\right\}\subset F_d$ and for $\zeta \in C$
\begin{equation}
\sup_{f\in F_d}\int_E \ln{f(x) \over \prob f(x)} \zeta(dx)\geq \max_{i=1,\ldots,p}\int_E\ln{f_i(x) \over \prob f_i(x)} \zeta(dx)> \inf_{\zeta \in C} I_d(\zeta) - \epsilon
\end{equation}
and since $\epsilon$ can be  arbitrarily small we obtain that
\begin{equation}
\lim_{d\to \infty} \inf_{\zeta \in C} I_d(\zeta) = \inf_{\zeta\in C}I(\zeta),
\end{equation}
which together with \eqref{aa10} completes the proof in the case when $C$ is compact.
When $C$ is closed in ${\cal P}(E)$ we follow the arguments of the proof of Theorem 4.4 in \cite{DV3}.
\endproof
We recall now Lemma 2.5 of \cite{DV1}
\begin{lemma}\label{linv}
For $\nu\in {\cal P}(E)$ we have that $I(\nu)=0$ if and only if $\nu$ is an invariant measure for the transition operator $\prob(x, \cdot)$.
\end{lemma}

We have
\begin{lemma}\label{alem3}
When there is a unique invariant measure $\mu$ for the transition operator $\prob(x, \cdot)$ and $\mu \notin C$, where $C\subset {\cal P}(E)$ is closed  and either $E$ is compact and or $E$ is a complete  separable metric space and \eqref{D.1} is satisfied then
\begin{equation}
\inf_{\zeta \in C} I(\zeta)>0.
\end{equation}
\end{lemma}
\proof Assume first that $E$ is compact. Then ${\cal P}(E)$ is also compact.
Notice that ${\cal P}(E)\ni \zeta \to I(\zeta)$ is l.s.c., and therefore in the case of compact $E$, when $I(\zeta_n)\to 0$ there is a subsequence $n_k$  and $\zeta\in {\cal P}(E)$ such that $\zeta_{n_k} \Rightarrow \zeta$ and $I(\zeta)=0$. Therefore by Lemma \ref{linv} measure $\zeta$ coincides with $\mu$, a contradiction. When $E$ is a complete separable metric space then by Lemma \ref{alem1} for each $q>0$ the set $C_q=\left\{\zeta\in {\cal P}(E): I(\zeta)\leq q\right\}$ is compact and the remaining part of the proof follows as in the case of compact space $E$.
\endproof

\bmhead{Acknowledgments} Research supported by Polish National Science Centre grant no.  2020/37/B/HS4/00120.


\end{document}